\documentclass{article}
\usepackage[english]{babel}

\usepackage{PRIMEarxiv}

\usepackage[utf8]{inputenc} 
\usepackage[T1]{fontenc}    
\usepackage{hyperref}       
\usepackage{url}            
\usepackage{booktabs}       
\usepackage{amsfonts}       
\usepackage{nicefrac}       
\usepackage{microtype}      
\usepackage{lipsum}
\usepackage{fancyhdr}       
\usepackage{graphicx}       
\graphicspath{{media/}}     

\usepackage{centernot}
\usepackage{mathtools}
\usepackage{ stmaryrd }
\usepackage{xcolor}

\usepackage{subfig}
\usepackage{graphicx}			
\usepackage{microtype} 			
\usepackage{ifthen}		    	
\usepackage[final]{pdfpages}    
\usepackage{lipsum}				
\usepackage{csquotes}
\usepackage{verbatim}

 \usepackage{tabularray}
 \UseTblrLibrary{diagbox}
\usepackage{diagbox}

\usepackage{verbatim}

\definecolor{bg}{rgb}{0.95,0.95,0.95}

\usepackage{cite}

\title{Conservation of the passband signal amplitude using a filter based on the Fast Fourier Transform algorithm}

\author{
  Flávio Dalossa Freire\\
  PPGMNE \\
  UFPR \\
  Curitiba\\
  \texttt{dflavio152@gmail.com} \\
   \And
  Isabel Gebauer Soares \\
  PGMEC \\
  UFPR \\
  Curitiba\\
  \texttt{isaa.isinhaa@gmail.com} \\
}

\begin{document}
\maketitle

\begin{abstract}
In this work, we propose an algorithm for a filter based on the Fast Fourier Transform (FFT), which, due to its characteristics, allows for an efficient computational implementation, ease of use, and minimizes amplitude variation in the filtered signal. The algorithm was implemented using the programming languages Python, R, and MATLAB. Initial results led to the conclusion that there was less amplitude loss in the filtered signal compared to the FIR filter. Future work may address a more rigorous methodology and comparative assessment of computational cost.
\end{abstract}


\section{Introduction}

Digital filtering is an essential technique in signal processing, used to remove noise, extract relevant information, and improve data quality in various applications, from digital communication to biomedical signal analysis. Digital filters, such as low-pass, high-pass, band-pass, and band-stop filters, play a crucial role in signal manipulation to meet specific frequency requirements. However, despite their wide applications and advantages, digital filters face inherent challenges that affect the accuracy and quality of the filtered signal. One of the most notable problems is amplitude loss, even if minimal, which can occur over time. This loss can compromise the signal's integrity, leading to a reduction in the fidelity of the processed signal. Additionally, digital filters may not be perfect in rejecting unwanted frequencies, allowing noise or interference components to pass through, which can further degrade the signal of interest.

The need to address these limitations is fundamental to advancing the effectiveness of digital filters. Improvements in filter design techniques, as well as the development of adaptive algorithms, are necessary to minimize amplitude loss and enhance unwanted frequency rejection.

This article proposes an algorithm to mitigate the problem of amplitude loss in the passband of the digitized signal to be filtered.

\section{Main Types of Filters according to frequency interval}
\label{sec:headings}

The main types of digital filters include:

\begin{enumerate}
    \item Low-Pass: \\
         Allows the passage of low frequencies and attenuates high frequencies. It is widely used in various applications, such as signal smoothing and high-frequency noise removal. However, the low-pass filter faces significant challenges, such as amplitude loss at cutoff frequencies and the inability to completely reject unwanted frequencies close to the cutoff frequency.
    \item High-Pass: \\
         Allows the passage of high frequencies and attenuates low frequencies. Used in applications that require the removal of low-frequency components, such as in audio signals to eliminate background noise.
    \item Band-Pass: \\
        Allows the passage of a specific range of frequencies and attenuates frequencies outside this range. Essential in communication systems to isolate signals of interest in an occupied frequency spectrum.
    \item Band-Stop: \\
        Attenuates a specific range of frequencies, allowing the passage of frequencies outside this range. Used to eliminate interference or noise at specific frequencies.
\end{enumerate}

\noindent Among these, the low-pass filter is often considered the foundation, as other types can be developed as variations of it. However, despite its importance and wide application, the low-pass filter is not without problems. One of the main challenges is the \textbf{amplitude loss} at frequencies close to the cutoff frequency, which can affect the accuracy of the processed signal.
Additionally, the filter's phase response can introduce distortions in the signal, particularly in time-critical signals.
Another significant problem is the \textbf{inaccuracy in the attenuation of unwanted frequencies}, where noise or interference components may still be present after filtering, compromising signal quality.

\section{The Algorithm}

The result of digitizing the signal is a data set, a vector that contains the respective magnitudes collected over time. Thus, such a numerical vector can be subjected to the Fast Fourier Transform (FFT) algorithm \cite{shin2008fundamentals}.
The basis of the algorithm proposed in this article is the FFT and its inverse, the Inverse Fast Fourier Transform (IFFT) algorithm. These algorithms, FFT and IFFT, are implemented in most of the programming languages commonly used in data processing and scientific computing. Here, we will cover Python, R, and MATLAB for implementing the proposed algorithm.

The algorithm described here has interesting properties, it is very simple, addresses the problem of amplitude loss in the passband of the filtered signal, and is computationally efficient as it is based on the FFT. The algorithm consists of three steps:

\begin{enumerate}
    \item Apply the FFT to the data set of the digitized signal (moving to the frequency domain).
    \item Remove the undesired frequencies or set the frequencies that should remain in the data array transformed to the frequency domain, skipping its first position.
    \item Apply the IFFT to the updated array, moving back to the time domain.
\end{enumerate}

In step 1, it is important to take the maximum number of sample points (signal measurements in time) to complete an integer number of cycles, thus mitigating the leakage problem inherent to the FFT. The leakage problem is minimized in large sample sizes in high-frequency signals, even with non-integer numbers of cycles, but it becomes more important in smaller samples and lower frequencies.

In step 2, it is crucial not to use mirroring or flipping methods to zero out the undesired frequencies in the transformed data vector to ensure that there are no other alterations to the original signal beyond the elimination of those undesired frequencies.

It is also important to note that step 2 can be used to introduce any values in any positions of the transformed vector. In the examples presented here, zeros will be introduced in the positions of the frequencies to be filtered out.

\begin{figure}[!htb]
    \centering
    \includegraphics[width=.825\textwidth]{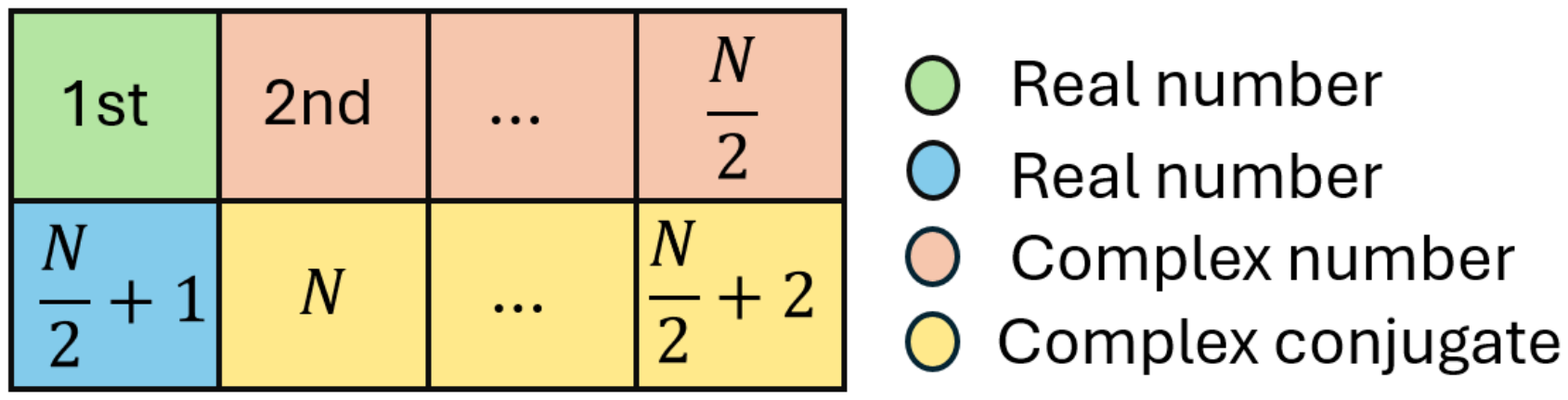}
    \caption{FFT Result structure for N even}
    \label{fig: n_even}
\end{figure}

\begin{figure}[!htb]
    \centering
    \includegraphics[width=.80\textwidth]{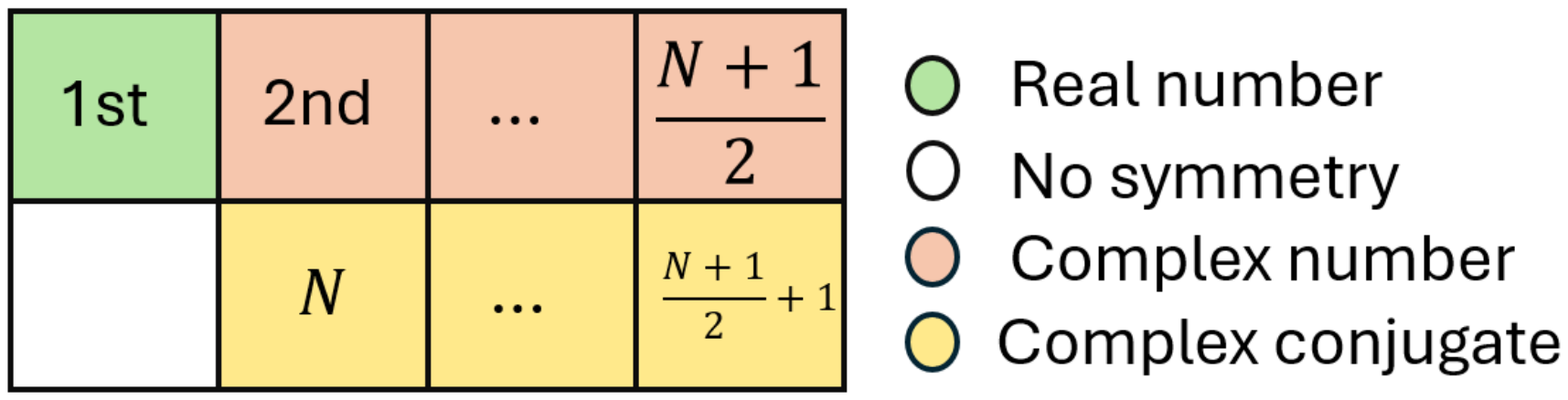}
    \caption{FFT Result structure for N odd}
    \label{fig: n_odd}
\end{figure}

The FFT algorithm returns a vector of the same size as the original data. For an even number of elements, N, this vector returned by the FFT will have a structure as schematized in \autoref{fig: n_even}, where the first element and the $N/2+1$ element are real numbers, and the rest are complex numbers. Among these complex numbers, there is a symmetry of conjugate pairs. The second element is conjugate with the last one, the third with the second last, and so on. For an odd number of elements, N, this vector returned by the FFT will have a structure as schematized in \autoref{fig: n_odd}. In this case, the first element is a real number and the rest are complex numbers. Among these complex numbers, there is a symmetry of conjugate pairs similar to the even N case. The second element is conjugate with the last one, the third with the second last, and so on.

\newpage

Each element of the vector transformed to the frequency domain has a correspondence in the time domain. The information contained in these elements can even be manipulated to obtain different results in the response vector returned by the inverse Fourier transform (IFFT). For example, a manipulation of the first element of the vector returned by the FFT results in a change in the intercept (the mean) of the vector reconverted to the time domain. A multiplication by two of the imaginary components in the already filtered vector, returned by the FFT, results in double the amplitude of the wave represented by the corresponding vector transformed by the IFFT. Thus, we can conclude that the first element of the vector returned by the FFT correlates with the mean of the data in the vector resulting from the IFFT, and information about the wave amplitude in the time domain is contained in the imaginary components of the fast transformed vector.

Thus, in step 2, the first element should always be disregarded to symmetrically zero out the frequencies in the transformed signal, starting from the second element. This is done by skipping the first element of the transformed vector and zeroing out the others whose indices correspond to the unwanted frequencies. These indices are obtained simply by taking the positions of the unwanted frequencies, or the extremes of their interval, in the corresponding frequency vector.  In this way, we have a \textbf{shifted symmetry} starting from the second element. If the element in the second position is to be zeroed, its conjugate pair in the last position must also be zeroed. Similarly, if the element in the third position is assigned a value of zero, its conjugate pair in the second last position must also be zeroed, and so on. Therefore, the algorithm works by assigning values to element $i$ and to its complex conjugate $N-i$, considering $i$ as the element index, $N$ as the total number of elements, and $0$ as the index of the first position.

We present examples of implementing the algorithm in the programming languages Python, R \cite{Rrr}, and Matlab \cite{MATLAB2024}. We also provide a suggestion in pseudo-code. Since our sample of signal points is centered at zero, in these codes we zero out the first element without worrying about changes in the intercept. However, in signals not centered at zero, this element should be preserved or restored to its original value before step 3 of the algorithm. Additionally, the suggestions presented here are not aimed at maximum computational efficiency but at illustrating the implementation of the algorithm with the best possible readability.

In this article, all comparisons and applications of the proposed algorithm, the option of \textbf{shifted symmetry} was used, argument \textit{shift\_sim}. The implemented functions have the option not to use shifted symmetry, considering the reader’s interest in comparing results using the example code.

\section{Pseudo-Code suggestion}

\begin{verbatim}
Function filter_IF(point_freq, interval, y, fs, shift_sim):
    
    If interval is None AND point_freq is None:
        Raise ValueError "At least one of interval or punctual frequencies parameters 
        must be provided"
    
    Apply FFT to original data points and store in Fft_smooth
    Copy Fft_smooth to Fft_smooth_p

    Set N to the length of y

    Calculate frequency vector far using N and fs
    
    far = (fs / N) * np.arange(N)  # Python exemple

    Get indices where far is less than or equal to interval[1] and store in indices1
    Get indices where far is greater than or equal to interval[0] and store in indices2
    Set index1 to the last element of indices1
    Set index2 to the first element of indices2

    If shift_sim is False (symmetric filtering):
        Set elements of Fft_smooth between (index1+1) and (N-index1-1) to 0
        Set elements of Fft_smooth up to index2 to 0
        Set elements of Fft_smooth from (N-index2) to the end to 0

    If shift_sim is True (shifted symmetric filtering):
        Set elements of Fft_smooth between (index1+1) and (N-index1) to 0
        Set elements of Fft_smooth up to index2 to 0
        Set elements of Fft_smooth from (N-index2+1) to the end to 0

    If point_freq is not None (point filtering):
        
        Find the index where far equals point_freq
        Set elements of Fft_smooth_p to 0 
        except for the index corresponding to point_freq and 
        its shifted symmetric position N - index

    Get the real part of the inverse FFT of Fft_smooth and store in y_cut
    Get the real part of the inverse FFT of Fft_smooth_p and store in y_cut_p

    Create a DataFrame df with columns:
        'FFt_filtered_i': Fft_smooth
        'FFt_filtered_p': Fft_smooth_p
        'y_cut': y_cut
        'y_cut_p': y_cut_p

    Return df
    
\end{verbatim}

\section{Code Examples}

Here, we present examples of the implementation of the aforementioned algorithm in the Python, R, and MATLAB programming languages.

\subsection{Python Example}
\begin{verbatim}
# Function filter_IF
def filter_IF(point_freq=None, interval=None, y=None, fs=None, shift_sim=False):

    if interval is None and poin_freq is None:
      
      raise ValueError("At least one of interval or punctual frequencies parameters 
      must be provided")    

    # Apply FFT in original data points   
    Fft_smooth = fft(y) 
    Fft_smooth_p = np.copy(Fft_smooth)

    # Number of points
    N = len(y)

    # Frequency vector
    far = (fs / N) * np.arange(N)  # Hz

    # Setting indexes
    indices1 = np.where(far <= interval[1])[0]
    indices2 = np.where(far >= interval[0])[0]
    index1 = indices1[-1]
    index2 = indices2[0]

    if not shift_sim: # Symmetric filtering
        
        Fft_smooth[(index1+1):(N-index1-1)] = 0
        Fft_smooth[:index2] = 0
        Fft_smooth[(N-index2):] = 0

    if shift_sim: # Shifted symmetric filtering
       
        Fft_smooth[(index1+1):(N-index1)] = 0
        Fft_smooth[:index2] = 0
        Fft_smooth[(N-index2+1):] = 0

    if point_freq is not None: # Point filtering
        
        index = np.where(far == point_freq)[0][0]
        Fft_smooth_p[~np.isin(np.arange(N), [index, N-index])] = 0    

    y_cut = np.real(ifft(Fft_smooth))  # Get the real part
    y_cut_p = np.real(ifft(Fft_smooth_p))  # Get the real part

    df = pd.DataFrame({
        'FFt_filtered_i': Fft_smooth,
        'FFt_filtered_p': Fft_smooth_p,
        'y_cut': y_cut,
        'y_cut_p': y_cut_p
    })

    return df
\end{verbatim}

\subsection{R Example}
\begin{verbatim}
## Function filter_IF 

filter_IF <- function(point_freq=NULL, interval=NULL, y, fs, shift_sim=FALSE) {
  
  Fft_smooth <- fft(y)
  Fft_smooth_p <- Fft_smooth 
  N <- length(y)
  
  # Frequency vector
  
  far <- (fs / N) * (0:(N-1))  # Hz
  
  if(!is.null(interval) & !shift_sim){
    
    indices1 <- which(far <= interval[2])
    indices2 <- which(far >= interval[1])
    index1 <- indices1[length(indices1)]
    index2 <- indices2[1]
    Fft_smooth[(index1+1):(N-index1)] <- 0
    Fft_smooth[1:index2] <- 0
    Fft_smooth[(N-index2+1):N] <- 0
    
  }
  if(!is.null(interval) & shift_sim){
    
    indices1 <- which(far <= interval[2])
    indices2 <- which(far >= interval[1])
    index1 <- indices1[length(indices1)]
    index2 <- indices2[1]
    Fft_smooth[(index1+1):(N-index1+1)] <- 0
    Fft_smooth[1:(index2-1)] <- 0
    Fft_smooth[(N-index2+3):N] <- 0
    
  }
  
  if (!is.null(point_freq)){
    
    index <- which(far == point_freq)
    Fft_smooth_p[-c(index,(N-index+2))] <- 0
    
  }
  
  else{
    
    message("At least one of interval or punctual frenquencies parameters 
            must be specified")
    
  }
  
  y_cut <- Re(fft(Fft_smooth, inverse = TRUE))/N # Get the real part
  y_cut_p <- Re(fft(Fft_smooth_p, inverse = TRUE))/N # Get the real part
  
  df <- data.frame(FFt_filtered_i=Fft_smooth,
                   FFt_filtered_p=Fft_smooth_p,
                   y_cut=y_cut, y_cut_p=y_cut_p)
  
  return(df)
}
\end{verbatim}

\subsection{MATLAB Example}
\begin{verbatim}
%% Function filter_IF

function df = filter_IF(point_freq, interval, y, fs, shift_sim)
    Fft_smooth = fft(y);
    Fft_smooth_p = Fft_smooth;
    N = length(y);
    
    % Frequency vector
    far = (fs / N) * (0:(N-1));
    
    if ~isempty(interval) && ~shift_sim

        indices1 = find(far <= interval(2));
        indices2 = find(far >= interval(1));

        index1 = indices1(end);
        index2 = indices2(1);

        Fft_smooth((index1+1):(N-index1)) = 0;
        Fft_smooth(1:index2) = 0;
        Fft_smooth((N-index2+1):N) = 0;

    end


    if ~isempty(interval) && shift_sim

        indices1 = find(far <= interval(2));
        indices2 = find(far >= interval(1));

        index1 = indices1(end);
        index2 = indices2(1);

        Fft_smooth((index1+1):(N-index1+1)) = 0;
        Fft_smooth(1:(index2-1)) = 0;
        Fft_smooth((N-index2+3):N) = 0;
    end

    
    if ~isempty(point_freq)
        index = find(far == point_freq);
        Fft_smooth_p(setdiff(1:N, [index, N-index+2])) = 0;
    else
        error(['At least one of interval or punctual' ...
            ' frequencies parameters must be specified']);
    end

    
    y_cut = real(ifft(Fft_smooth)); % Get the real part
    y_cut_p = real(ifft(Fft_smooth_p)); % Get the real part
    
    df.FFt_filtered_i = Fft_smooth;
    df.FFt_filtered_p = Fft_smooth_p;
    df.y_cut = y_cut;
    df.y_cut_p = y_cut_p;
end
\end{verbatim}


\section{Discussion and Results}
\label{sec:others}

\subsection{Signal Amplitude: Filtered versus Theoretical}

\begin{figure}[!htb]
    \centering
    \includegraphics[width=.51\textwidth]{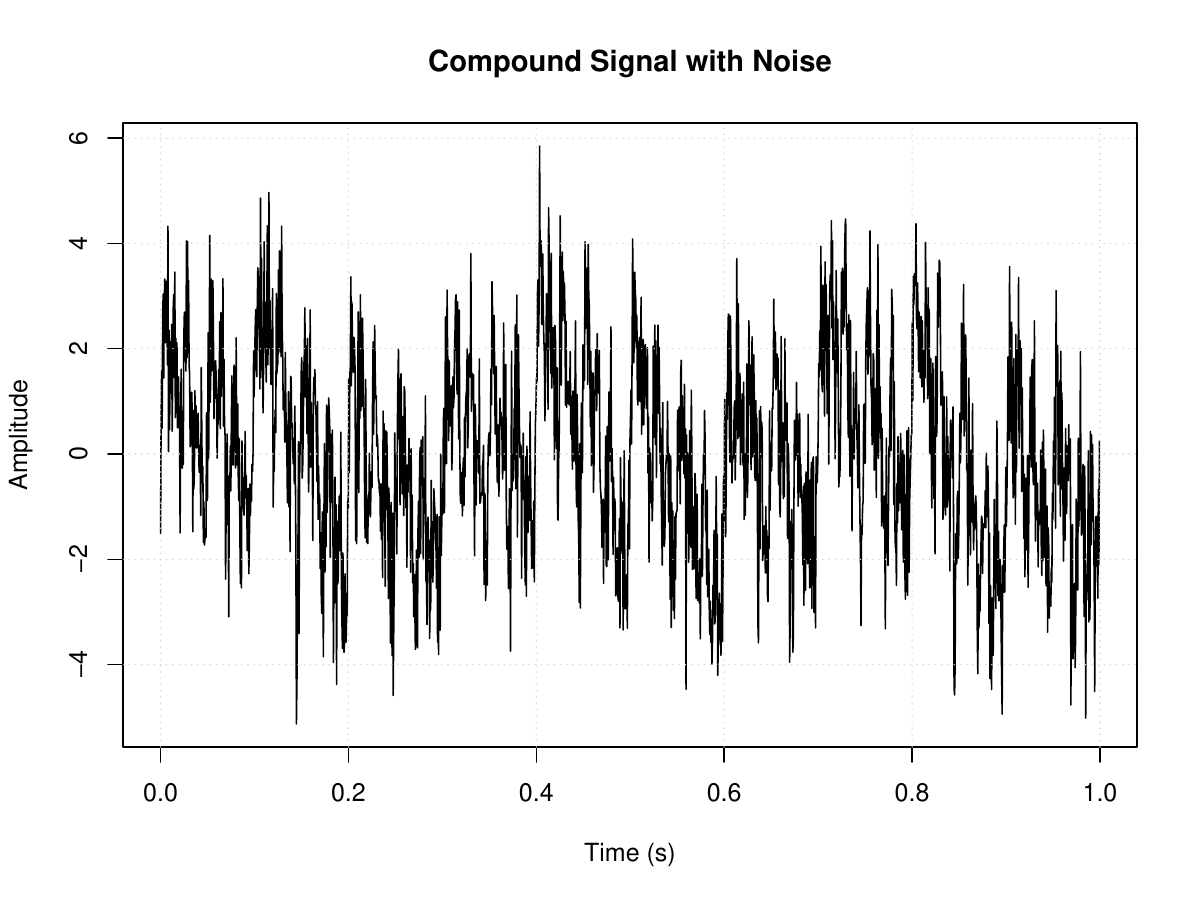}
    \caption{Signal with Noise Visualization}
    \label{fig: plot0}
\end{figure}

The R function presented above was applied to a data set representing a signal generated by simulation, summing 5 sine waves with frequencies 3, 10, 20, 40, and 80 Hz, amplitudes 1, with a sampling frequency of 1770 Hz over an interval of 1 second, resulting in 1777 simulated measurements. White noise with a normal distribution, mean 0, and standard deviation 1 was added to these measurements, \autoref{fig: plot0}.

Next, the amplitudes of the signals filtered by frequency intervals and by specific values were compared with the theoretical components of the simulated signal, and the amplitude variation in each case was observed. Furthermore, the results for band-pass filtering and single-frequency filtering were compared.

\begin{figure}[!htb]
    \centering
    \includegraphics[width=.8\textwidth]{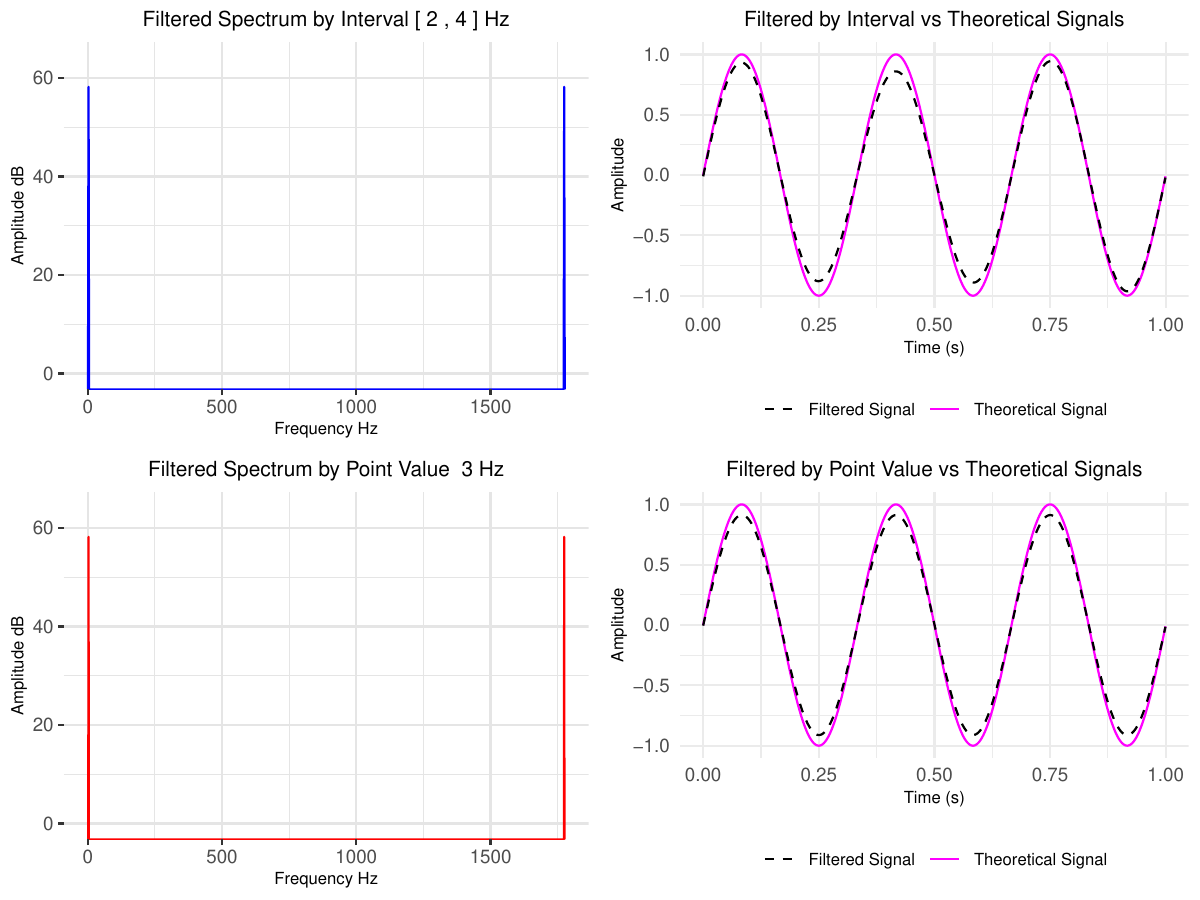}
    \caption{Comparison between band-pass filter $2$ to $4$ Hz and single-pass $3$ Hz}
    \label{fig: plot1}
\end{figure}

\begin{figure}[!htb]
    \centering
    \includegraphics[width=.8\textwidth]{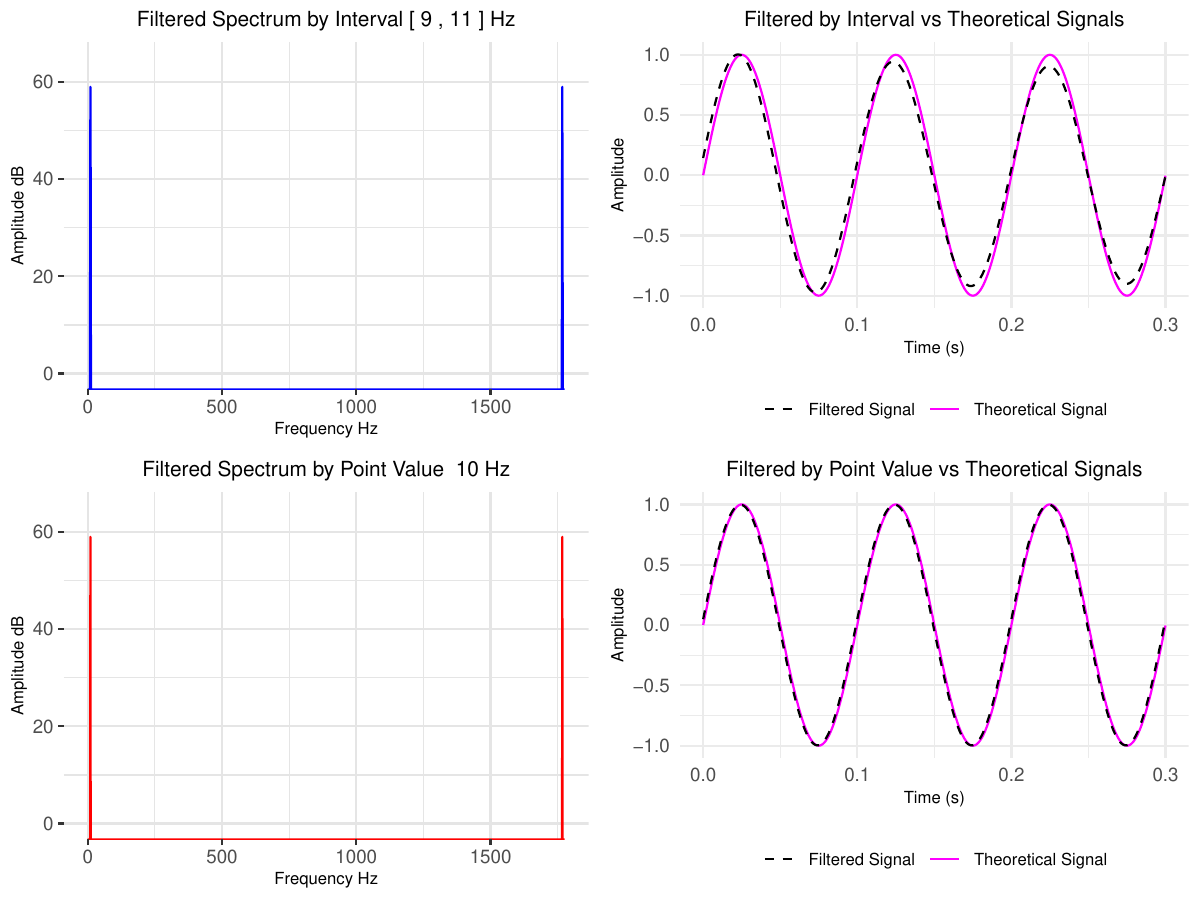}
    \caption{Comparison between band-pass filter $9$ to $11$ Hz and single-pass $10$ Hz}
    \label{fig: plot2}
\end{figure}

\begin{figure}[!htb]
    \centering
    \includegraphics[width=.8\textwidth]{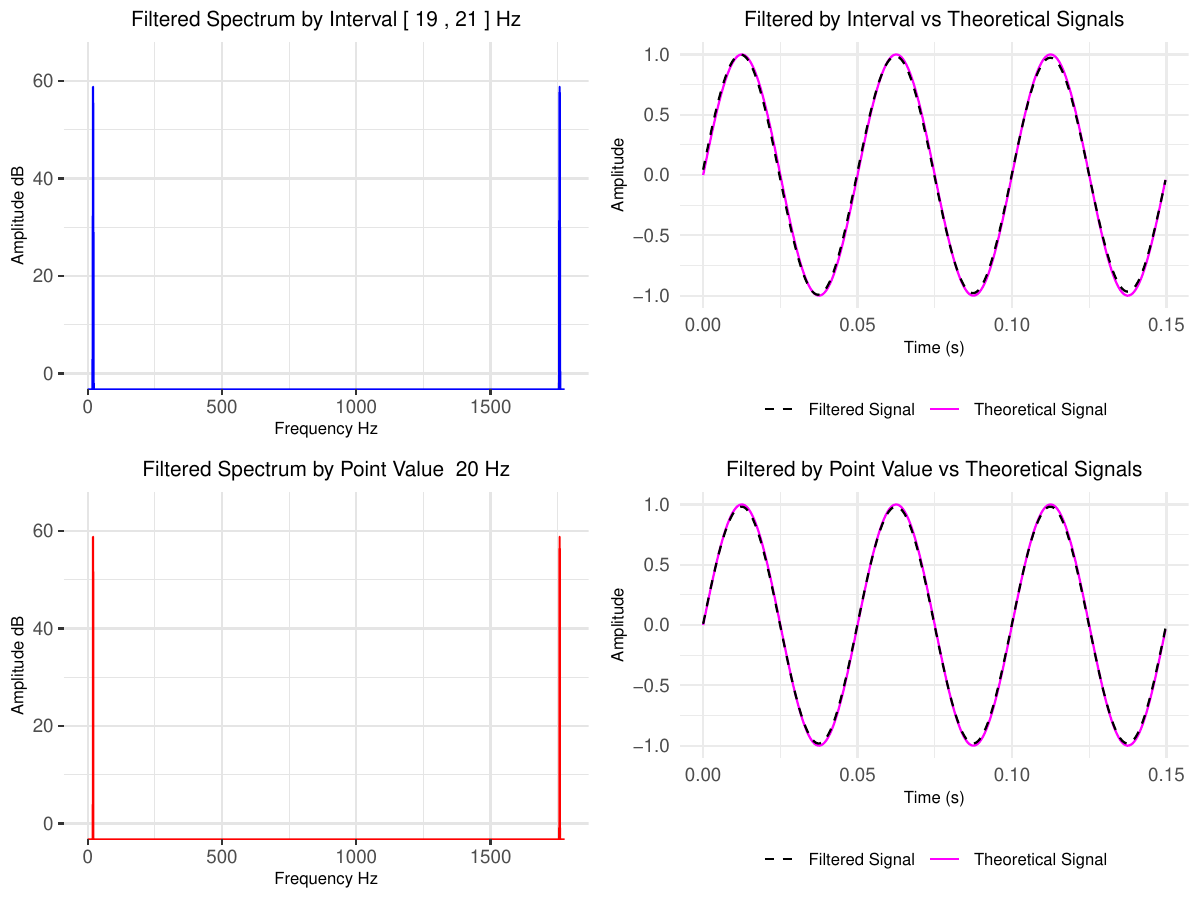}
    \caption{Comparison between band-pass filter $19$ to $21$ Hz and single-pass $20$ Hz}
    \label{fig: plot3}
\end{figure}

\begin{figure}[!htb]
    \centering
    \includegraphics[width=.83\textwidth]{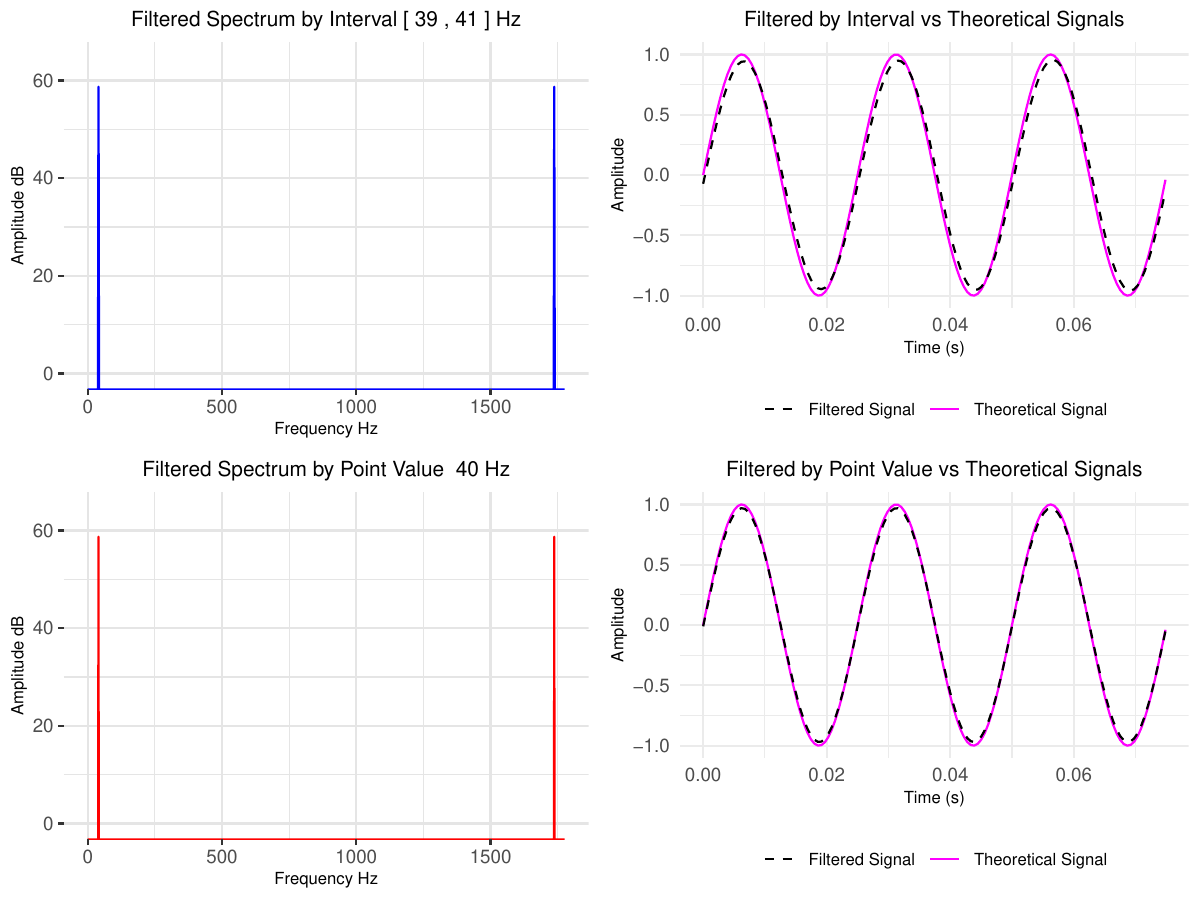}
    \caption{Comparison between band-pass filter $39$ to $41$ Hz and single-pass $40$ Hz}
    \label{fig: plot4}
\end{figure}

\begin{figure}[!htb]
    \centering
    \includegraphics[width=.85\textwidth]{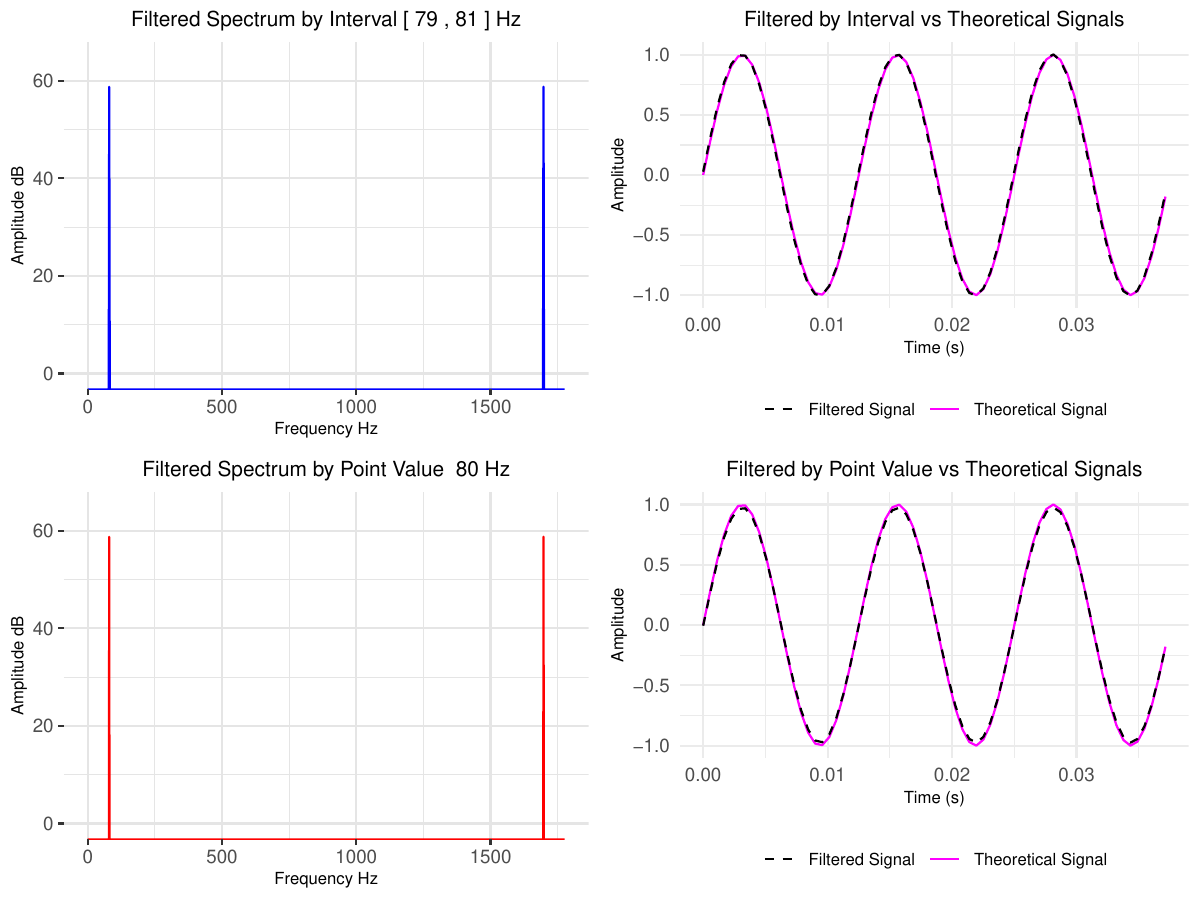}
    \caption{Comparison between band-pass filter $79$ to $81$ Hz and single-pass $80$ Hz}
    \label{fig: plot5}
\end{figure}

A slight decrease in the amplitude of the filtered signal at 3 Hz was observed in the case of this simulated signal with its respective noise \autoref{fig: plot1}. This variation tends to decrease with an increase in the passband frequency, as shown in Figures \ref{fig: plot1} to \ref{fig: plot5}.

\newpage

\subsection{Signal Amplitude: FFT Filter versus FIR Equiripple Filter}

\begin{figure}[!htb]
    \centering
    \includegraphics[width=.7\textwidth]{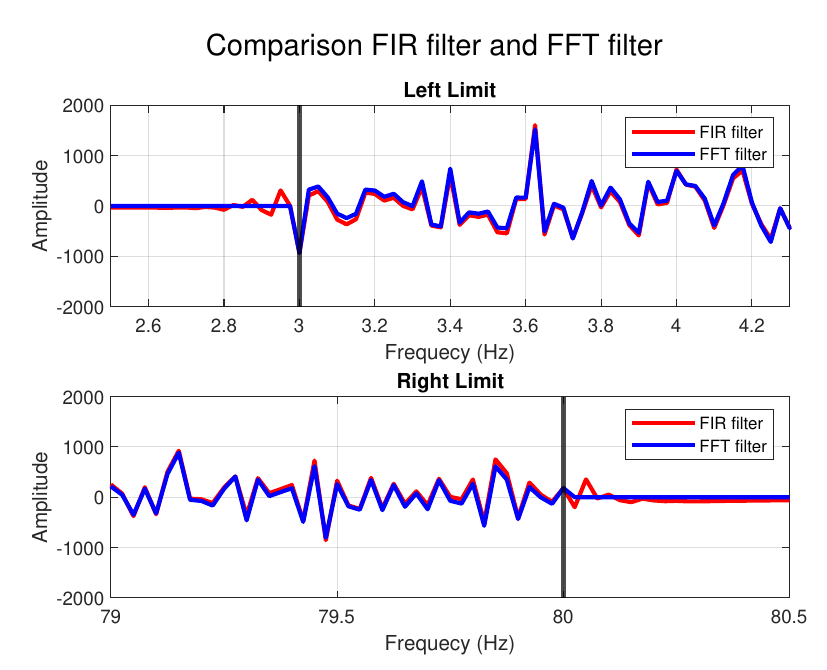}
    \caption{Comparison between FFT filter and FIR filter for the 3 to 80 Hz band in the frequency domain}
    \label{fig: comparacao3a80}
\end{figure}

\begin{figure}[!htb]
    \centering
    \includegraphics[width=.8\textwidth]{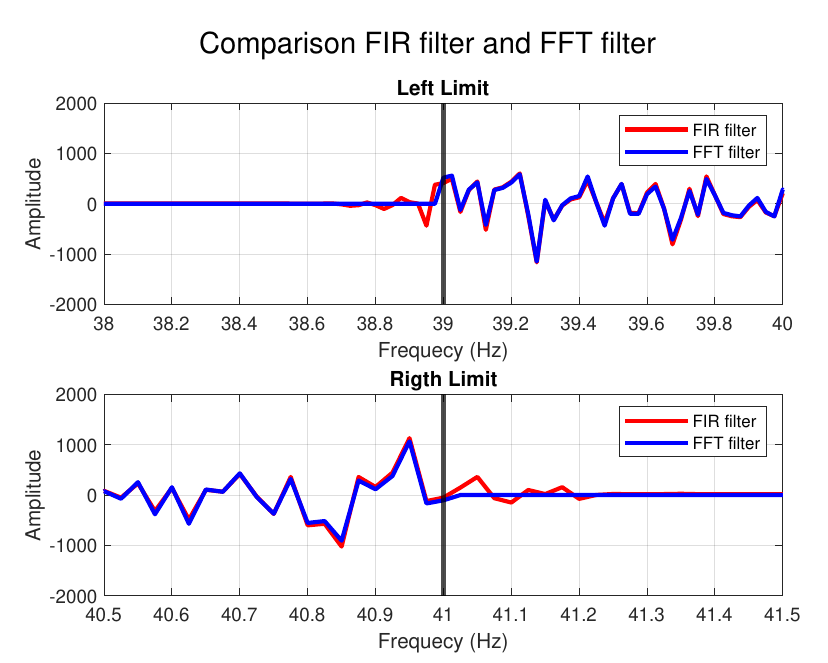}
    \caption{Comparison between FFT filter and FIR filter for the 39 to 41 Hz band in the frequency domain}
    \label{fig: comparacao39a41}
\end{figure}

\begin{figure}[!htb]
    \centering
    \includegraphics[width=.8\textwidth]{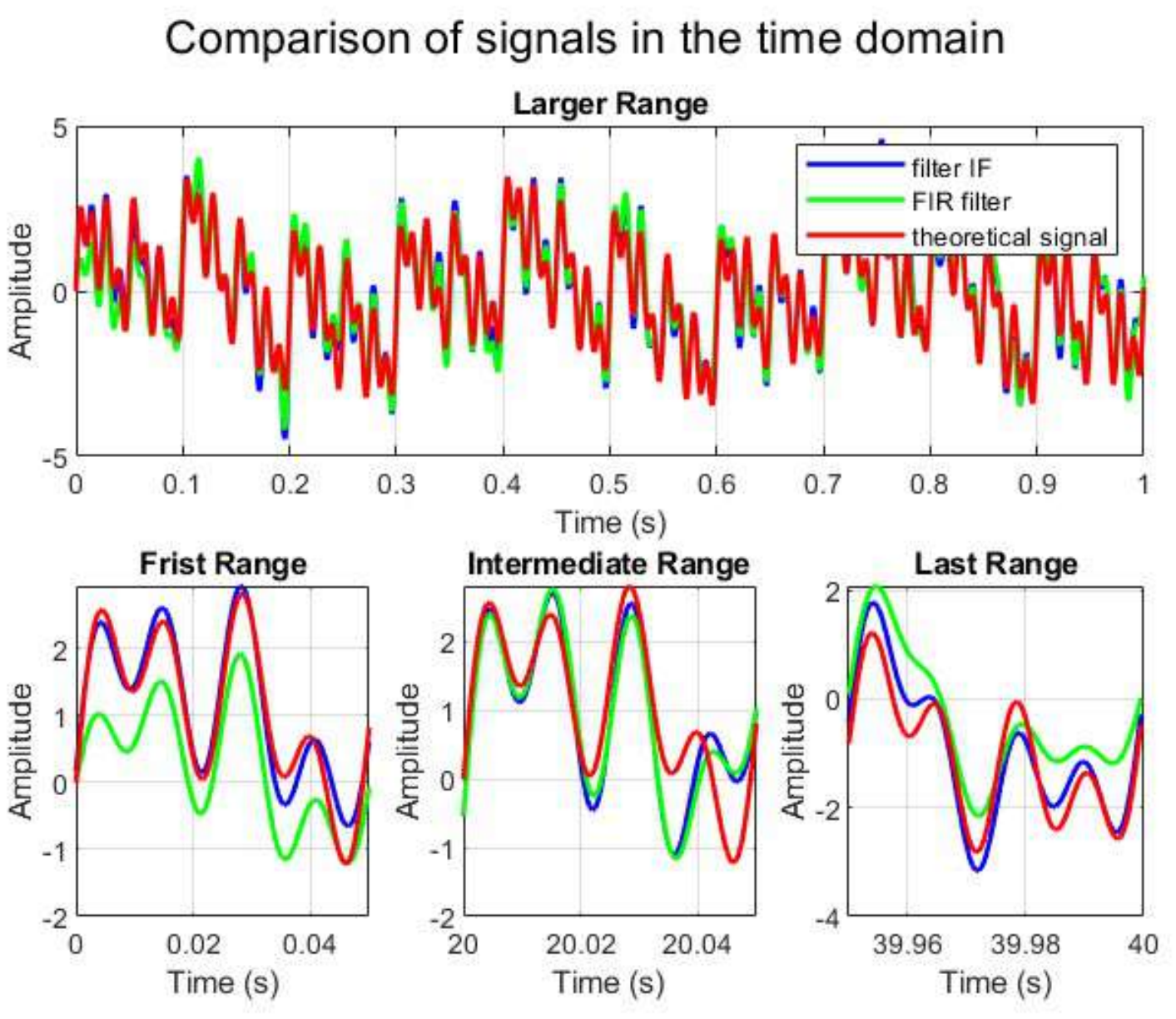}
    \caption{Comparison between FFT filter and FIR filter for the 3 to 80 Hz band in the time domain}
    \label{fig: comparacaoTempo1}
\end{figure}

\begin{figure}[!htb]
    \centering
    \includegraphics[width=.79\textwidth]{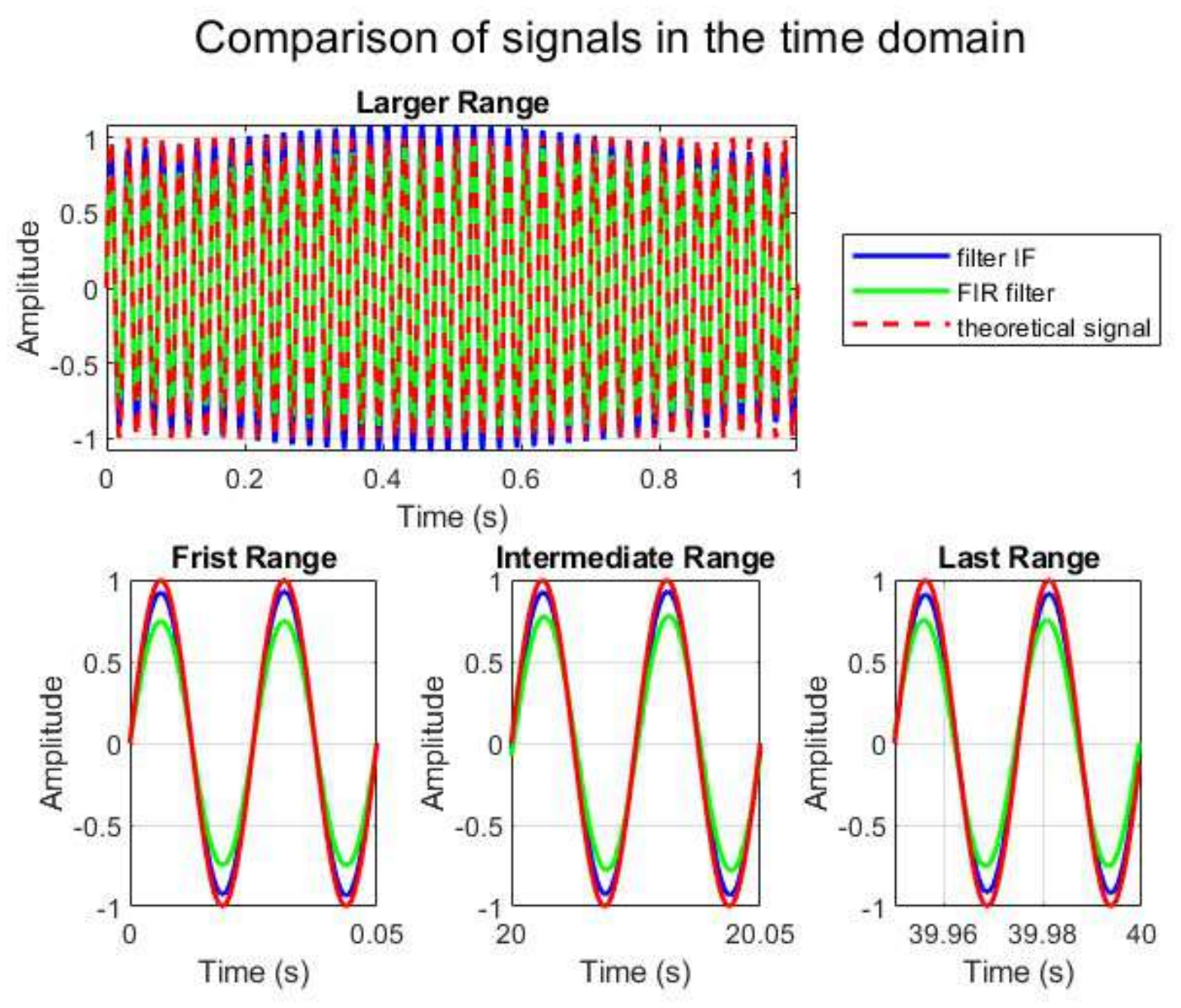}
    \caption{Comparison between FFT filter and FIR filter for the 39 to 41 Hz band in the time domain}
    \label{fig: comparacaoTempo2}
\end{figure}

Finally, the implementation in MATLAB of the algorithm-based filter proposed here (FFT Filter) was compared with a design using an equiripple FIR filter (FIR Filter). The passbands chosen for both filters were 3 to 80 Hz and 39 to 41 Hz. Thus, the accuracy of frequency domain cutoffs (\autoref{fig: comparacao3a80}) and (\autoref{fig: comparacao39a41}), adherence to the theoretical signal, and amplitude variation in the time domain were observed for these two bands, 3 to 80 Hz and 39 to 41,  Hz (\autoref{fig: comparacaoTempo1}) and (\autoref{fig: comparacaoTempo2}).

\newpage

Below, in MATLAB code, we have the parameters of equiripple FIR band-pass filters, used for 3 to 80 Hz band and 39 to 41 Hz band. In \autoref{tabelaRMSE}, the results for the root mean square error (RMSE) for both filters, relative to the theoretical signal for each considered band, are shown.

\begin{itemize}
    \item MATLAB code for FIR filter, 3 to 80 Hz 
\end{itemize}

\begin{verbatim}

function Hd = FIRfilter1
%FIRfilter1 Returns a discrete-time filter object.

% MATLAB Code
% Generated by MATLAB(R) 23.2 and Signal Processing Toolbox 23.2.
% Generated on: 13-Jul-2024 17:39:59

% Equiripple Bandpass filter designed using the FIRPM function.

% All frequency values are in Hz.
Fs = 2000;  % Sampling Frequency

Fstop1 = 2.5;             % First Stopband Frequency
Fpass1 = 3;               % First Passband Frequency
Fpass2 = 80;              % Second Passband Frequency
Fstop2 = 80.5;            % Second Stopband Frequency
Dstop1 = 0.001;           % First Stopband Attenuation
Dpass  = 0.057501127785;  % Passband Ripple
Dstop2 = 0.0001;          % Second Stopband Attenuation
dens   = 20;              % Density Factor

% Calculate the order from the parameters using FIRPMORD.
[N, Fo, Ao, W] = firpmord([Fstop1 Fpass1 Fpass2 Fstop2]/(Fs/2), [0 1 ...
                          0], [Dstop1 Dpass Dstop2]);

% Calculate the coefficients using the FIRPM function.
b  = firpm(N, Fo, Ao, W, {dens});
Hd = dfilt.dffir(b);

% [EOF]

\end{verbatim}

\begin{itemize}
    \item MATLAB code for FIR filter, 39 to 41 Hz 
\end{itemize}

\begin{verbatim}

function Hd = FIRfilter2
%FIRfilter2 Returns a discrete-time filter object.

% MATLAB Code
% Generated by MATLAB(R) 23.2 and Signal Processing Toolbox 23.2.
% Generated on: 13-Jul-2024 18:28:18

% Equiripple Bandpass filter designed using the FIRPM function.

% All frequency values are in Hz.
Fs = 2000;  % Sampling Frequency

Fstop1 = 38.5;            % First Stopband Frequency
Fpass1 = 39;              % First Passband Frequency
Fpass2 = 41;              % Second Passband Frequency
Fstop2 = 41.5;            % Second Stopband Frequency
Dstop1 = 0.001;           % First Stopband Attenuation
Dpass  = 0.057501127785;  % Passband Ripple
Dstop2 = 0.001;           % Second Stopband Attenuation
dens   = 20;              % Density Factor

% Calculate the order from the parameters using FIRPMORD.
[N, Fo, Ao, W] = firpmord([Fstop1 Fpass1 Fpass2 Fstop2]/(Fs/2), [0 1 ...
                          0], [Dstop1 Dpass Dstop2]);

% Calculate the coefficients using the FIRPM function.
b  = firpm(N, Fo, Ao, W, {dens});
Hd = dfilt.dffir(b);

% [EOF]

\end{verbatim}

\begin{table}[!htb]
\centering
\begin{tabular}{llll} 
\hline
\multicolumn{1}{c}{\diagbox{{}measure}{signal}} & FIR filter & FFT filter & Unfiltered Signal  \\ 
\hline
RMSE: 3 to 80 Hz                                      & 163.0046   & 155.5202   & 566.9388           \\
RSME 39 to 41 HZ                                      & 39.3361    & 24.1285    & 693.1644           \\
\hline
\end{tabular}
 \caption{RSME Comparison between filters and unfiltered signal as reference}
 \label{tabelaRMSE}
\end{table}

One of the reasons for the better performance of the FFT filter using shifted symmetry can be seen in \autoref{fig: comparacao3a80} and \autoref{fig: comparacao39a41}. The proposed algorithm allows for a perfectly rectangular window, with a straight and vertical cut at the chosen limits for the signal in the frequency domain. Note that there are no residual frequencies beyond the cut-off limits of the FFT filter using shifted symmetry. Thus, there is no need for a trade-off between the precision in eliminating unwanted frequencies and maintaining the amplitude of the filtered signal.

\section{Conclusion}

It is observed that the RMSE value, used as a metric for comparing the FIR Filter and the FFT Filter, was lower in the case of the FFT filter implemented by the algorithm proposed here. The amplitude obtained in the 39 to 41 Hz band was also higher for the signal processed by the FFT filter, as seen from inspection of the graph in (\autoref{fig: comparacaoTempo2}). This finding, combined with a better RMSE value, already allows us to conclude that the amplitude loss was less for the FFT filter. Nevertheless, for a more rigorous demonstration in future work, it is necessary to use another metric better suited to measure amplitude variation across a larger number of data simulations, especially varying the generation of white noise, both in comparison with the theoretical signal and between filters. Therefore, future work may address a more rigorous methodology and comparative assessment of computational cost, including different scenarios and a more detailed description of the theory behind the proposed algorithm.

\bibliographystyle{unsrt}  
\bibliography{references}

\end{document}